% ******************** D E F I N I T I O N S ****************

\magnification=1200
\parskip = 6 pt

\def\real{I\kern-4pt R}

\def\ol{{\overline{l}}}
\def\nnor{{\overline{r}}}
\def\ogamma{{\overline{\gamma}}}
\def\ob{{\overline{b}}}
\def\verysmall{{\vskip 0.22 cm }}
\def\ms{{\medskip }}
\def\bb{{\cal B}}

\def\ss{{\cal S}} 
\def\aa{{\cal A}}
\def\rr{{\cal R}}
\def\oP{{\overline{P}}}
\def\of{{\overline{f}}}

\def\cmn{{C_m\times C_n}}
\def\cr{{\hbox{cr}}}
\def\;{{ \ }}

\def\onm{\hbox{\hglue 0.03 cm }{\ominus_{\hbox{\hglue -0.06 cm } {}_{n} }\hbox{\hglue -0.01 cm}}}
\def\onp{\hbox{\hglue 0.03 cm }{\oplus_{\hbox{\hglue -0.06 cm } {}_{n} }\hbox{\hglue -0.01 cm}}}
\def\okm{\hbox{\hglue 0.03 cm }{\ominus_{\hbox{\hglue -0.06 cm } {}_{k} }\hbox{\hglue -0.01 cm}}}
\def\okp{\hbox{\hglue 0.03 cm }{\oplus_{\hbox{\hglue -0.06 cm } {}_{k} }\hbox{\hglue -0.01 cm}}}
\def\oom{{\ominus_{\hbox{\hglue -0.06 cm } {}_{m} }\hbox{\hglue -0.05 cm}}}
\def\oop{{\oplus_{\hbox{\hglue -0.06 cm } {}_{m} }\hbox{\hglue -0.05 cm}}}

\def\xx{{\cal X}}
\def\tt{{\cal T}}
\def\cc{{\cal C}}
\def\yy{{\cal Y}}
\def\dd{{\cal D}}
\def\ii{{\cal I}}
\def\mybox#1{{\ \hbox{#1}\ }}
\def\subsecc#1{{\verysmall\noindent{\it #1}\verysmall}}

\def\boxit#1{\vbox{\hrule\hbox{\vrule\kern3pt
      \vbox{\kern3pt#1\kern2pt}\kern2pt\vrule}\hrule}}
\setbox4=\vbox{\hsize 1 pt }

\def\hal{\vrule height 5 pt width 0.05 in  depth 0.8 pt}
\font\boluno = cmbx10   scaled    \magstep 1

% ***********************************************************
% ********************  F O N T S ***************************

 0

\font\bolcer = cmbx10   scaled    \magstep 0

\font\bitcer = cmbxti10 scaled    \magstep 0
\font\tty    = cmtt10 scaled \magstep 0

 0
\font\bituno = cmbxti10 scaled    \magstep 1

\font\scacer = cmcsc10 scaled \magstep 0
 0

\font\scacerotro = cmssbx10 scaled \magstep 0
\font\scacerotro = cmfib10 scaled \magstep 0
\font\scacerotro = cmssbx10 scaled \magstep 0
\font\ssscer     = cmss10 scaled \magstep 0

\footline={}

% ***********************  T I T L E ************************

\vglue 1.0 cm

\centerline{\boluno The conjecture 
{\bituno cr}({\bituno C}${}_{\hbox{\bitcer m}} 
                    \hskip 0.1 cm \times \hskip -0.1 cm 
                    $ {\bituno C}${}_{\hbox{\bitcer n}}$)
	{\boluno =} ({\bituno m} -- 2){\bituno n} \hglue 0.1 cm   is true}
\smallskip
\centerline{{\boluno for all but finitely many 
 {\bituno n}, for each {\bituno m}}}

\vskip 1.0 cm

\centerline{\scacer Lev Yu.~Glebsky}
\vskip 0.4 cm
\centerline{\it Department of Mechanics and Mathematics,}
\centerline{\it
		Nizhny Novgorod State University,
		Russia}
		
\bigskip
\centerline{and}
\bigskip

\centerline{\scacer Gelasio Salazar\footnote{*}{\rm 
Corresponding Author. E--mail: {\tty gsalazar@cactus.iico.uaslp.mx} }}
\vskip 0.4 cm
\centerline{\it IICO--UASLP, San Luis Potosi, Mexico}

\vskip 1.0 cm

\centerline{15 September 2000}

% ************************************************************

% ********************  A B S T R A C T **********************

\vglue 2.5 cm 

{\leftskip = 20 pt \rightskip = 20 pt

\noindent {\bolcer Abstract.}  
It has been long conjectured that 
the crossing number of $\cmn$ 
is $(m-2)n$, for all $m, n$ such that $n \ge m \ge 3$.  
In this paper it is proved that this conjecture holds 
for all but finitely many $n$, for each $m$.  
More specifically, it is shown that if 
$n \ge (m/2)((m+3)^2/2 + 1)$ and $m \ge 3$, then 
the crossing number of $\cmn$ is exactly
$(m-2)n$, as conjectured. The 
proof is largely based on the theory of arrangements, 
introduced by Adamsson and further
developed by Adamsson and Richter. 
\par}

\vglue 2 cm 

\centerline{\ssscer To be  submitted to the Journal of Graph Theory}

% ************************************************************

\vfill\eject

 \pageno=1

 \footline={\centerline{\the\pageno}}

% ***************** S E C T I O N     1 *********************

% ************** I n t r o d u c t i o n *********************

\vskip 1.0 cm

\centerline{\scacerotro 1. INTRODUCTION}

\vskip 0.35 cm

% ***** I n t r o d u c t o r y      P a r a g r a p h *****

	In 1973, Harary, Kainen, and Schwenk proved
	that toroidal graphs can have arbitrarily large
	crossing numbers [7].  In the same paper, they
	put forward the following conjecture.

\proclaim Conjecture [HKS--Conjecture]. 
	The crossing number $\cr(\cmn)$ of the Cartesian product
	$\cmn$ is $(m-2)n$, for all $m, n$
	such that $n \ge m \ge 3$. 

	This has been proved for $m, n$ 
	satisfying $n \ge m$, $m \le 7$ [12, 5, 6, 
	11, 9,  3, 10, 4, 1].  Our aim in this paper is to 
	show that the HKS--conjecture holds for all but finitely
	many $n$, for each fixed $m \ge 3$.

\proclaim Main Theorem.
	Let $m, n$ be integers such that
	$n \ge 
	(m/2)((m+3)^2/2 + 1)$, $m \ge 3$. Then
	$\cr(\cmn) = (m-2)n$.

        Although we do not use the notion of arrangement 
	explicitly, the proof of the Main Theorem is largely based on 
	the theory of arrangements, introduced by 
	Adamsson 
	[1], and further developed by Adamsson and Richter [2].

	In the proof we make frequent use of the Jordan Curve
	Theorem.  With this exception, the proof is
	self--contained. 
	As we point out in the last section, 
	the statement of the Main Theorem can be slightly
	improved using the 
	general bound $\cr(\cmn)\ge (m-2)n/2$ [8].

	The heart of the proof of the Main Theorem is 
	the following.

\proclaim Theorem 1.
	Let $m, n$ be integers such that $n \ge m \ge 3$.
	Then every robust drawing of
	$\cmn$ has at least
	$(m-2)n$ crossings.  
	
	Roughly speaking (formal definitions are 
	in Section 2), a drawing of $\cmn$ is robust 
	if (i) for every three $m$--cycles $R, R', R''$, there 
	is a component of $\real^2\setminus R$ that intersects 
	both $R'$ and $R''$; and
	(ii) to every $m$--cycle $R$ we can assign two 
	disjoint $m$--cycles $R', R''$, both disjoint from $R$, 
	such that every cycle between (with respect to
	a circular relation defined below) $R'$ and $R$
	is disjoint from $R''$.

	The argument that shows that the Main Theorem 
	follows from Theorem 1 can be outlined as
	follows. Let $m \ge 3$ be fixed, and let
	$n_0 = (m+3)^2/2 + 1$.  It is easy to check 
	that the Main Theorem 
	is a consequence of  the following auxiliary 
	statement: for all 
	$n \ge n_0$, $\cr(\cmn) \ge  \min\{(m-2)n, m(n-n_0)\}$.  
	This statement is proved by induction on $n$.
	The base case is $n = n_0$, for which there is
	nothing to prove.  Suppose the statement is
	true for $n = k-1 \ge n_0$, and let $\dd$ be a
	drawing of $C_m\times C_k$.  Since $k$ is large
	(enough) compared to $m$, then either $\dd$ is 
	robust or some $m$--cycle has $m$ or more 
	crossings.
	In the first case the statement follows directly
	(without using the inductive assumption) from
	Theorem 1.  In the second case, the statement follows
	by  applying the
	inductive assumption to the drawing obtained
	by removing $R$ from $\dd$.

	The complete, formal proof that the Main Theorem 
	is a consequence of Theorem 1 is given in Section 8.
	The rest of 
	this paper is devoted to the proof of Theorem 1.

	The strategy of the proof of Theorem 1 is
	to show that in every
	robust drawing of $\cmn$, 
	we can associate to each of the $n$ $m$--cycles 
	at least $m-2$ crossings, in such a way that no
	crossing is associated to more than one 	
	$m$--cycle.   
	
	In Section 2 we introduce basic definitions, notation,
	and terminology.
	In Section 3 we analyze drawings 
	of structures that consist of three closed curves
	plus a set of arcs that meet 
	the curves in the
	same order.   The major topological results 
	needed in the proof of Theorem 1 (namely 
	Corollaries 3 and 5) are established in this	
	section; the rest of the proof consists mostly 
	of combinatorial arguments.
	We chose to present these topological
	results at this early stage in order to prevent
	a disruption of the discussion in later sections.

	In Section 4 we prove some
	basic facts on robust drawings. In Section 5 we 
	specify the set of crossings associated to
	each $m$--cycle in a robust drawing.  In Section 6
	we show that no crossing is associated to more than
	one $m$--cycle, and in Section 7 we show that 
	there are at least $m-2$ crossings associated to
	each $m$--cycle.  Finally, in Section 8 we prove
	Theorem 1 and the Main Theorem.  Section 9 contains
	some final remarks.

% ***************** S E C T I O N     2 **********************

% ***********  Basic definitions and Terminology   ***********

\vskip 0.8 cm

\centerline{\scacerotro 2. BASIC DEFINITIONS, NOTATION, AND TERMINOLOGY}

\vskip 0.35 cm

% ***** I n t r o d u c t o r y      P a r a g r a p h *****

\subsecc{2.1 \  Addition, subtraction, and
		circular relation on $Z_k$}

	For each integer $k > 1$, 
	we denote addition and subtraction on $Z_k$
	by the symbols $\okp$ and $\okm$, respectively.
	We define the {\it circular relation} $\preceq$ on 
	$Z_k$ by the rule $i \preceq j$ iff 
	$j \okm i \in \{0,\ldots, \lfloor{k/2}\rfloor\}$. 
	We write $i \prec j$ if $i \preceq j$ and
	$i \neq j$.   This relation has the following properties:
	(i) if $i \not\preceq j$, then $j \prec i$ (we remark
	that if $k$ is even, then it is possible that 
	$i \prec j$ and $j \prec i$ ); 
	(ii) if $i \prec j \prec i \okp l$ for some $l \neq 0$,
	then $j = i \okp x$ for some 
	$x\in \{1,\ldots,l\okm 1\}$; and 	
	(iii) if $i, j < n/2$, then 
	$i \preceq j$ iff $i \le j$.

\subsecc{2.2 \ The Cartesian product $\cmn$}

	The {\it Cartesian product} $C_m \times C_n$ is a 
$4$--regular graph with $mn$ vertices $v(i,j)$, where 
$0 \le i \le m - 1$ and $0 \le j \le n - 1$. The 
vertices are labeled so that
the vertices adjacent to
$v(i,j)$ are $v(i\oom 1,j), v(i\oop 1,j), v(i,j\onm 1),$ and 
$v(i,j\onp 1)$.

The edge set of $C_m \times C_n$ is naturally partitioned 
into $m$ edge sets of $n$--cycles and $n$ edge sets of 
$m$--cycles. To help comprehension,
we color the $n$--cycles blue and the $m$--cycles 
red. We label the blue cycles $(v(i,j)), j\in Z_n$, by
$B(i)$, $i\in Z_m$, and the red cycles $(v(i,j))$, 
$i\in Z_m$, by $R(j)$, $j\in Z_n$.

Let $i\in Z_m, j,k \in Z_n, j\neq k$. 
The blue edge that joins $v(i,j\onm 1)$ to $v(i,j)$
is denoted $bl(i,j)$. 
The {\it open blue path} $P(i,j,k)$ is 
the sequence of edges and vertices
$bl(i,j\onp 1), v(i,$ $j\onp 1),
 bl(i,j\onp 2), \ldots, 
v(i,k\onm 1), bl(i,k)$.
The {\it closed blue path} $\oP(i,j,k)$ is obtained
by adding $v(i,j)$ at the beginning and $v(i,k)$ at
the end of $P(i,j,k)$.

\subsecc{2.3 \ Arcs, well--behaved collections of arcs,
	tangential intersections, crossings}

	An {\it open arc} $\gamma$ is the image of a 
	local homeomorphism $f:(0,1) \to \real^2$ (a 
	{\it building} map for
	$\gamma$) with the property that the unique continuous 
	extension $\of$ of $f$ to $[0,1]$ is such that
	$\of(0)\neq\of(1)$ and $\of(0), \of(1)\notin\gamma$.
	Denote by $\ogamma$
	the image of $\of$.  The points $\of(0)$ and $\of(1)$
	are the {\it end points} of both $\gamma$ and $\ogamma$. 
	A {\it closed arc} $\gamma$ is the image of a 
	local homeomorphism
	$f:S^1\to \real^2$ (a {\it building} map for
	$\gamma$).  If $\gamma$ is an (open or closed)
	arc that has some one--to--one building map, then
	$\gamma$ is {\it simple}.

	Let $\gamma$ be an open arc with building map
	$f:(0,1)\to \real^2$.  An arc $\delta$ is a 
	{\it subarc} of $\gamma$ if there are $a,b$,
	$0\le a < b \le 1$, such that the map
	$g:(0,1)\to \real^2$ defined by the rule
	$g(x) = f(x(b-a)+a)$ is a building map for $\delta$.
	If $0< a < b < 1$, then $\delta$ is a 
	{\it totally proper subarc} (or simply {\it tp--subarc}) 
	of $\gamma$.  Thus, if $\delta$ is a
	$tp$--subarc of $\gamma$, then
	no endpoint of $\gamma$ is an endpoint of $\delta$.

	Let $\cc$ be a collection of closed arcs, and let 
	$C, C', C''\in \cc$.  If no component of
	$\real^2\setminus C$ intersects both $C'$ and
	$C''$, then $C$ {\it separates} $C'$ {\it from}
	$C''$.  If no arc in $\cc$ separates two
	arcs in $\cc$ from each other, then 
	$\cc$ is {\it nonseparating}.

	Let $\gamma$ be an arc
	with building map $f:X\to \real^2$.  
	For each $z\in\gamma$, the
	{\it multiplicity} of $z$ is
	$|\{ y \in X \ | \ f(y) = z \}|$.  
	If $z\in \gamma$ has multiplicity $\mu > 1$, then
	$z$ is a {\it self--intersection of $\gamma$ of
	multiplicity $\mu$}.  It is easy to see that 
	the multiplicity of a self--intersection of  $\gamma$
	is independent of the building map chosen for $\gamma$. 
	Clearly, $\gamma$ is simple iff it has no 
	self--intersections.

% 	*****************************************************

	A collection $\aa$ of arcs is
	{\it well--behaved} if, for every $\gamma, \delta$ in $\aa$:

\item{(i)} 
	every self--intersection of $\gamma$ has multiplicity $2$,
	and has a neighborhood that contains no other
	self--intersections of $\gamma$;
	\vglue -0.18 cm 
\item{(ii)} $\gamma$ and $\delta$ intersect each other 
	a finite number of times; and 	\vglue -0.22 cm 
\item{(iii)}
	every intersection between $\gamma$ and $\delta$
	is a self--intersection 
	of neither $\gamma$ nor $\delta$.

	Let $\{\gamma,\delta\}$ be a well--behaved collection
	of arcs, and let $x$ be an 
	intersection point between $\gamma$ and $\delta$
	(a self--intersection if $\gamma = \delta$).
	Then there is a set $N\subset \real^2$ homeomorphic to a 
	closed disc, whose boundary $\partial{N}$ is a simple 
	closed arc,
	and whose interior $N^o$ contains $x$, such that
	$(\gamma\cup\delta)\cap N^o$ is the union of
	two open 
	arcs $\alpha, \beta$, with the following properties:
	(i) $\alpha$ and $\beta$ are simple $tp$--subarcs
	of $\gamma$ and $\delta$, respectively;
	(ii) $\overline{\alpha}\cap\overline{\beta}=\{x\}$;
	and 
	(iii) $(\gamma\cup\delta)\cap {N} = 
	\overline{\alpha}\cup\overline{\beta}$. 
	If $x$ can be removed by an isotopy on
	$\overline{\alpha}$ totally contained 
	in ${N}$, that leaves $\partial{N}$ fixed, 
	then $x$ is {\it tangential}. Otherwise $x$  is 
	a {\it crossing}. 

	Suppose that $\gamma=\delta$ (so $x$ is a 
	self--intersection), and $x$ is a crossing.
	Let $a, a'$ and $b, b'$ be the end points of $\alpha$ 	
	and $\beta$, respectively.  There are 
	(uniquely determined) simple arcs
	$\phi, \psi$ contained in 
	$\partial{N}$  such that
	(i) the end points of both $\phi$ and $\psi$ 
	are in $\{a, a', b, b'\}$; 
	(ii) $\overline{\phi}\cap\overline{\psi}=
	\emptyset$;
	and
	(iii) 	$\gamma'= 
	(\gamma\setminus(\alpha\cup\beta))\cup(\phi\cup\psi)$
	is an arc. This arc $\gamma'$ has exactly one
	fewer self--intersection than $\gamma$.  We say
	that $\gamma'$ is obtained from $\gamma$ by
	{\it smoothing out} the self--intersection $x$.

% 	*****************************************************

\subsecc{2.4 \ Drawings of graphs, definition of \ $\sqcap_\dd$}

	A {\it drawing} $\dd$ of a simple graph $G$ is a
	representation of $G$ in the plane such that: 
	(i) each vertex is represented by a point, and 
	no two different vertices are represented by the 
	same point; (ii) each edge $e$ is 
	represented by an open arc,  so
	that the end points of the representation of
	$e$ are precisely the points that represent
	the vertices incident with $e$; 
	(iii) no representation of an edge contains
	a representation of a vertex.

\noindent{\bf Remark.} 
	For simplicity, if there is only one drawing under	
	consideration, 
	we often make no distinction between
	a substructure of the graph (such as a vertex, or
	a path, or a cycle)  and the subset of $\real^2$
	that represents it.  Throughout this work, we
	have taken special care to ensure that no confusion
	arises from this practice.

	Let $\dd$ be a drawing of $\cmn$.
	Suppose that each of $H$ and $K$ is either 
	an open path or a red cycle, and that no edge is in
	both $H$ and $K$.
	Denote by 	$H\sqcap_{\dd} K$
	the set of pairwise intersections of edges in $\dd$ that 
	involve
	one edge in $H$ and one edge in $K$.
	If $\dd$ is the only drawing under
	consideration, we omit the reference to $\dd$ and
	simply write
	$\sqcap$.

	A drawing of a simple graph is {\it good} if (i) no
	edge has a self--intersection; (ii) no two
	adjacent edges intersect; (iii) no two
	edges intersect each other
	more than once; and (iv) each intersection of
	edges is a crossing.

	The {\it crossing number} $\cr(G)$ of a graph $G$ is 
	the minimum number of pairwise intersections of 
	edges in a drawing of $G$ in the plane.  An
	{\it optimal} drawing of $G$ is a drawing where
	the crossing number is attained.  It is
	a routine exercise to show that every
	optimal drawing of a graph is a good drawing	
	(hence the term {\it crossing} number, in view of (iv)).

\subsecc{2.5 \ Robust drawings of $\cmn$}

	Let $\dd$ be a drawing of $\cmn$. Fix $j\in Z_n$.
If $R(j)\sqcap_\dd R(k) = \emptyset$ for some red cycle $R(k)$, 
then let
$
b(\dd,j) = \min \{b\in Z_{n}\ |\ R(j \onm  b)\sqcap_\dd R(j)=\emptyset 
	\}. 
$
If $b(\dd,j)$ is defined, and there is a red cycle 
$R(l)\notin\{R(j\onm b(\dd,j)), R(j\onm ( b(\dd,j) - 1)),\ldots,
	R(j)\}$ such 
that 
$R(j \onm c) \sqcap_\dd R(l) =\emptyset$ for each
$c$ such that $0 \le c \le b(\dd,j)$, then define
	
$
a(\dd,j)=\min\{a\in Z_n \ | \ \forall  c, 0 \leq c \leq b(\dd,j), 
	R(j \onm c)\sqcap_\dd R(j\onp  a)=\emptyset\}.
$ 	
	
If there is only one drawing under consideration, we
omit the reference to $\dd$ and simply write 
$b(j)$ and $a(j)$. 

	A drawing $\dd$ of $\cmn$ is {\it red--nonseparating}
	if $\{R(0),\ldots,R(j-1)\}$ is nonseparating,
	and it is {\it relaxed} if 
		
	for every $j\in Z_n$,
	$a(\dd,j), b(\dd,j)$ are defined and
	$a(\dd,j) + b(\dd,j) < n/2$. 
	If $\dd$ is relaxed, then
	$\bb(\dd) = \max\{ a(\dd, j), b(\dd, j) \ | \ j\in Z_n\}$. 	
	Finally, $\dd$ is {\it robust} if it is
	red--nonseparating and relaxed.

% ***************** S E C T I O N     4 **********************

% ***************  Topological  Prerequisites    *************

\vskip 0.5 cm

\centerline{\scacerotro 3. ANALYSIS OF CROSSINGS IN 
	(3,s)--CONFIGURATIONS}

\vskip 0.35 cm

% ***** I n t r o d u c t o r y      P a r a g r a p h *****

	Let $k \ge 0$  and $s \ge 1$ be integers.
	Let $\cc = \{C_0,C_1,C_2\}$ be a collection 
	of closed arcs,
	and let $\aa = \{A_0,\ldots,A_{s-1}\}$ be a
	collection of  open arcs.  The
	pair $(\cc,\aa)$ is a
	{\it k--intersecting $(3,s)$--configuration} if
	the following are satisfied:

\ms

\item{(i)}
	$\cc\cup\aa$ is well--behaved;\vglue -0.15 cm 
\item{(ii)} $\cc$ is nonseparating;\vglue -0.15 cm 

\item{(iii)}
	$|C_0 \cap C_1| = k$, and
	$C_2$ is disjoint from $C_0\cup C_1$; \vglue -0.15 cm 
\item{(iv)}
	each $A_i$ has one end
	point (the {\it initial} point 
	$t(i)$ of $A_i$) in $C_0$, and
	the other end point  
	 (the {\it final} point $f(i)$ of $A_i$)
	in $C_2$;\vglue -0.15 cm 
\item{(v)} each $A_i$ intersects 
	$C_1$ in exactly one point 
	(the {\it middle} point $w(i)$ of $A_i$), and 
	does so tangentially; \vglue -0.15 cm 
\item{(vi)}
	with the exception of the intersection points mentioned
	in (v), all the intersections in $\cc\cup\aa$ 
	(including self--intersections) are crossings.

\ms

	For each $A_i$ in $\aa$, let $T_i$ denote
	the subarc of $A_i$ whose end points are 
	$t(i)$ and $w(i)$, and 
	let $F_i$ denote the subarc of $A_i$ whose
	end points are 
	$w(i)$ and 
	$f(i)$.  The subarcs $T_i$ and $F_i$ 
	are the {\it initial} and {\it final} subarcs
	of $A_i$, respectively.  Thus, $A_i =
	T_i \cup F_i \cup \{w(i)\}$.

	A $(3,s)$--configuration is {\it clean} if
	for every $A_i \in \aa$, (i) $T_i$ does not
	intersect  $C_2$,  and (ii) $F_i$ does not
	intersect $C_0$.  An intersection between 
	arcs in a 
	$(3,s)$--configuration is {\it good} if it occurs 
	either between a  $T_i$ and  an
	$F_j$, for some $i\neq j$, or between two 
	initial subarcs $T_i, T_j$, for some $i\neq j$.
	To emphasize the fact that 
	every good intersection is necessarily a
	crossing,  we also use the term {\it good crossing} 
	to refer to a good intersection.

	\proclaim Lemma 2. Every clean $0$--intersecting
	$(3,s)$--configuration 
	has at least $s - 2$ good crossings. 

	Let $(\cc,\aa) = (\{C_{0}, C_1, C_2\},\{A_0,\ldots,A_{s-1}\})$ 
	be a clean 
	$(3,s)$--configuration.	
	If $s < 3$, then there is nothing to prove. So 
	we assume $ s \ge 3$.  Suppose the statement is true
	for $s = 3$ (that is, every clean 
	$(3,3)$--configuration has
	at least one good crossing).  An elementary counting 
	argument then shows that if $s \ge 4$, then every 
	$(3,s)$--configuration has at least 
	${s\choose 3}/(s-2) \ge s - 2$ good crossings.  Thus 
	it suffices to show that every clean 
	$(3,3)$--configuration has
	at least one good crossing.
	Therefore we assume 
	$\aa= \{A_0, A_1, A_2\}$.  

	Each $T_i$ has a unique
	subarc $T_i'$ such that (i) one end point of $T_i'$ is
	in $C_0$; (ii) the other end point of $T_i'$ is $w(i)$;
	and (iii) $T_i'$ does not intersect $C_0$. 
	Similarly, 
	each $F_i$ has a unique
	subarc $F_i'$ such that (i) one end point of $F_i'$ is
	in $C_2$; (ii) the other end point of $F_i'$ is $w(i)$;
	and (iii) $F_i'$ does not intersect $C_2$. 
	For each $i$, let $A_i'= T_i' \cup F_i' \cup \{w(i)\}$, 
	and let $\aa' = \{A_0', A_1', A_{2}'\}$. 
	It is easy to check that $(\cc,\aa')$ is also a clean
	$(3,3)$--configuration.  Moreover, every good crossing
	of $(\cc,\aa')$ is a good crossing of
	$(\cc,\aa)$.  Thus it suffices to show that
	$(\cc,\aa')$ has at least one good crossing.

	By construction, no $T_i'$ intersects $C_0$. Since
	$(\cc,\aa')$ is clean, no
	$T_i'$ intersects $C_1\cup C_2$.  Thus, no
	$T_i'$ intersects $C_0 \cup C_1 \cup C_2$.  
	Similarly, no $F_i'$ intersects
	$C_0 \cup C_1 \cup C_2$.  

% ------------------------------------------------------------

	We claim we can assume that 
	$C_0, C_1$, and $C_2$  have no 
	self--intersections.  For suppose that
	$x$ is a self--intersection of $C_i$,
	for some $i\in\{0,1,2\}$.
	Then $x$ has multiplicity $2$, and $x$ is 
	the end point of no arc in $\aa$.  Thus 
	we can obtain a new closed arc $C_i'$ by
	smoothing out $x$ from $C_i$,
	without modifying any
	arc in $\cc\cup\aa$ other than $C_i$, so that
	$((\cc\setminus\{C_i\})\cup \{C_i'\} , \aa)$ is a clean 
	$(3,3)$--configuration 
	with the
	same number of good crossings as $(\cc,\aa)$.  

	Since an intersection involving $T_i'$ and $T_j'$ with
	$i\neq j$ is 
	good, we may assume $T_0', T_1', T_2'$ are pairwise disjoint.  
	Since $C_2$ does not 
	intersect $C_0\cup C_1\cup T_0'\cup T_1'\cup T_2'$, 
	it follows
	that $C_2$  is contained in one of the five components of 
	$\real^2\setminus(C_0\cup C_1\cup T_0'\cup 
	T_1'\cup T_2')$.  
	
	The boundary of one of these five
	components (say $U$) is $C_0$, and the boundary of 
	another of these
	five components (say $V$) is $C_1$.  Since no $F_i'$ intersects
	$C_0$, it follows that $C_2$ cannot be contained in $U$.
	Since each $w(i)$ is a tangential intersection, it
	follows that $C_2$ cannot be contained in $V$.
	Thus the boundary $\partial{W}$ of the component $W$ that
	contains $C_2$ is the disjoint union of 
	$T_i', T_j'$ (for some $i\neq j$), one 
	$tp$--subarc of $C_0$, one $tp$--subarc 
	of $C_1$, and $\{w(i),w(j), t(i),t(j)\}$.   
	Let $k$ be the integer in $\{0,1,2\}$ 
	different from $i$ and $j$. 
	Since $w(k)$ is not in $\partial{W}$, it follows that
	$F_k'$ must intersect (cross) $\partial{W}$.  Since
	$F_k'$ does not intersect $C_0\cup C_1$, then
	$F_k'$ must intersect $T_i'\cup T_j'$.  
	Since such an intersection is good, we are done. \hal\ms

% ------------------------------------------------------------

\proclaim Corollary 3.
	Let $(\cc,\aa) = (\{C_0,C_1,C_2\},\{A_0,\ldots,A_{s-1}\})$
	be a 
	$0$--intersecting $(3,s)$--confi\-gu\-ration. 
	Let $x_1$ denote the number of good crossings of
	$(\cc,\aa)$, let $x_2$ denote the number of
	initial arcs that cross $C_2$, and let $x_3$ 
	denote the number of final arcs that
	cross $C_0$.  Then 
	$x_1 + x_2 + x_3  \ge s - 2$.

\noindent{\bitcer Proof.} This follows from
	the definition of a clean configuration 
	and Lemma 2. \hal\ms

\proclaim Lemma 4.  
	Let $(\cc,\aa) = (\{C_0,C_1,C_2\},\{A_0,\ldots,
	A_{s-1}\})$ be a clean $k$--intersecting 
	$(3,s)$--configuration, with $k > 0$. 
	Then $(\cc,\aa)$ has at
	least $s - k$ good
	crossings.  

\noindent{\bitcer Proof.}  Using similar techniques as
	in the proof of Lemma 2, construct a  
	set $\aa'=\{A_0', A_1',$ $\ldots,$ $A_{s-1}'\}$ of open arcs
	such that (i) $(\cc,\aa')$ is a clean $k$--intersecting
	$(3,s)$--configuration; (ii) the initial arc
	$T_i'$ of each $A_i'$ does not intersect 
	$C_0\cup C_1\cup C_2$; 
	(iii) the final arc
	$F_i'$ of each $A_i'$ does not intersect 
	$C_0\cup C_1\cup C_2$; 
	and (iv) every good crossing of $(\cc,\aa')$ is
	a good crossing of $(\cc,\aa)$. 

	Thus, it suffices to show that $(\cc,\aa')$ has
	at least $s-k$ good crossings.
	
	As in the proof of Lemma 2, we can assume that
	$C_0, C_1,$ and $C_2$ have no self--intersections.
	Let $W$ denote the component of $\real^2\setminus
	(C_0\cup C_1)$ that contains $C_2$.  
	Let $\{D_1, D_2, \ldots, D_r\}$ denote the
	pairwise disjoint (necessarily  $tp$--) 
	subarcs of $C_0$ on $\partial{W}$ (clearly, 
	there is at least one such segment). 
	For each $j\in\{1,\ldots,r\}$, let $D_j'$ be 
	a  $tp$--subarc of $D_j$  
	that contains all the 
	$t(i)$'s in $D_j$.  
	Since $C_0$ is a simple arc, each $D_j'$ is 
	also a simple arc.

	For each $j\in\{1,\ldots,r\}$, let $\aa_j'$ 
	denote the set of
	arcs $A_i'$ in $\aa'$ such that $t(i) \in D_j'$. 
	By assumption, (i) no $T_i'$ and no $F_i'$ intersects
	$C_0\cup C_1\cup C_2$, and (ii) the only 
	intersection between each $A_i'$ and $C_1$ is 
	tangential.  It is easy to check that it follows
	that each $A_i'$ is 
	contained in $W$.  In particular, 
	each $A_i'$ has its initial point $t(i)$ in
	some $D_j'$.  Thus the collections $\aa_j'$
	partition $\aa'$.   Let $s_j = |\aa_j'|$.

	For each $j\in\{1,\ldots,r\}$, 
	let $U_j$ denote the unique component
	of $\real^2\setminus (C_0\cup C_1)$  different from 
	$W$ that contains $D_j'$ in its boundary. 
	For each $D_j'$, 
	draw an open arc $E_j$ very 
	close to $D_j'$, contained in $U_j$, 
	with the same endpoints as $D_j'$.
	Let $H_j$ denote the closed arc that consists
	of $D_j'\cup E_j$ plus the (common) endpoints of 
	$D_j'$ and $E_j$.  For each $j$ such that 
	$\aa_j'\neq\emptyset$, $(\{H_j, C_1, C_2\}, \aa_j')$ 
	is a $0$--intersecting $(3,s_j)$--configuration.  
	Each such $(3,s_j)$--configuration has at least
	$s_j - 2$ good crossings, by Lemma 2.
	Therefore
	$(\cc,\aa')$ has at least
	$\sum_{j=1}^r (s_j - 2)
	= -2r + \sum_{j=1}^r s_j$ good crossings.

	It is readily checked that the total number $k$ of 
	crossings between $C_0$ and $C_1$ is at least
	$2r$, and so 
	$-2r \ge -k$.  On the other hand, the collections
	$\aa_j'$ partition $\aa'$, and so
	$\sum_{j=1}^r s_j = \sum_{j=1}^r |\aa_j'| = |\aa'| 
	= s$.  Thus $(\cc,\aa')$ has at least
	$s - k$ good crossings. \hal\ms

\proclaim Corollary 5.
	Let $(\cc,\aa) = (\{C_0,C_1,C_2\},\{A_0,\ldots,A_{s-1}\})$
	be a 
	$k$--intersecting $(3,s)$--configuration, where $k > 0$. 
	Let $x_1$ denote the number of good crossings of
	$(\cc,\aa)$, let $x_2$ denote the number of
	initial arcs that cross $C_2$, and let $x_3$
	denote the number of final arcs that
	cross $C_0$.  Then
	$x_1 + x_2 + x_3 + k \ge s$. 

\noindent{\bitcer Proof.} This follows from the definition  
	of a clean configuration and Lemma 4. \hal\ms

% ***************** S E C T I O N     5 **********************

% ***********  Definition of the \ii_j             ***********

\vskip 0.5 cm

\centerline{\scacerotro 4. 
		PRELIMINARY RESULTS ON ROBUST DRAWINGS}

\vskip 0.35 cm
\ms

% ***** I n t r o d u c t o r y      P a r a g r a p h *****

\centerline{\it Throughout this section, 
	$\dd$ is a fixed robust good drawing of $\cmn$.}  

\ms

	We have the following preliminary observations	
	and conventions.
\ms

\item{(1)} As 
	$\dd$ is the only drawing under consideration,
	we shall omit the reference to $\dd$ in 
	the parameters $a(\dd,j), b(\dd,j)$, and $\bb(\dd)$
	(which are defined for all $j$, since $\dd$
	is robust), and in the symbol $\sqcap_\dd$,
	and simply write $a(j), b(j)$, $\bb$, and $\sqcap$,
	respectively.

\item{(2)}
	Since $\dd$ is good, all the intersections of edges
	are crossings.  To emphasize this, we do not speak
	of intersections of edges but of crossings of edges. 

\item{(3)} 
	If $H, K$ are (nonnecessarily different, nonnecessarily
	disjoint) subgraphs of $\cmn$, then we say that
	$H$ {\it crosses} $K$ if some edge of $H$ crosses some
	edge of $K$. 
\ms

	An $(i,j)$--crossing is a crossing between 
	an edge in $B(i)$ and an edge in $R(j)$. 
	If $B(i)$ crosses $R(j)$, then 
	the {\it first $(i,j)$--crossing from $v(i,j')$}
	(respectively {\it last})
	is the first (respectively last) $(i,j)$--crossing 
	we find as we traverse $B(i)$ 
	starting at $v(i,j')$ and finding the other
	vertices in $B(i)$ in the order
	$v(i,j'\onp 1), v(i,j'\onp 2),\ldots, v(i,j'\onm 1)$.

\noindent{\bf Remark.} For each  $R(j)$,
there is a unique component $\Omega_j$ of
$\real^2\setminus R(j)$ that intersects every red
cycle different from $R(j)$.  To see this, first
we note that, since $\dd$ is robust, then there is
an $R(k)$ such that $R(j)\sqcap R(k)=\emptyset$.
Then $\Omega_j$ is the component 
of $\real^2\setminus R(j)$ that contains
$R(k)$.  Indeed, if
$R(l)$ does not intersect $\Omega_j$
for some $l\neq j$, 
then $R(j)$ separates $R(k)$ from $R(l)$, contradicting
the assumption that $\dd$ is robust.

	For each $j\in Z_n$, let  
	$\Phi_j = \real^2\setminus (\Omega_j\cup R(j))$.

\noindent{\bf Remark.} 
 If $R(j)\sqcap R(k) = \emptyset$, then clearly
	$R(k)$ is contained in $\Omega_j$.

\proclaim Proposition 6. Let $j_1, j_2, 
	\in Z_n$, and let 
	$a_1, b_1, a_2, b_2\in Z_n$ satisfy 
	$0 < b_1 \le b(j_1), 0 < a_1 \le a(j_1)$, 
	$0 < b_2 \le b(j_2)$, and $0 < a_2 \le a(j_2)$. 
	Then the following hold for every $i\in Z_m$:\vglue -0.05 cm 
\item{(i)} If $j_2 = j_1 \onm b_1$, then  
	$P(i,j_1,j_1\onp a_1)$ and
	$P(i,j_2\onm \bb(\dd),j_2)$ 
	have no common edges. \vglue -0.2 cm 
\item{(ii)} If $j_2 = j_1 \onp a_1$, then  
	$P(i,j_1\onm b_1,j_1)$ and $P(i,j_2,j_2\onp \bb(\dd))$
	have no common edges. \vglue -0.2 cm 
\item{(iii)} If 
	$j_1 \onm b_1    = 	j_2 \onp a_2$, then
	$P(i,j_2\onm b_2,j_2)$ and $P(i,j_1,j_1\onp a_1)$
	have no common edges.\vglue -0.2 cm 
\item{(iv)} If $j_2 \preceq j_1 \onm b_1$ and $j_2\prec j_1$, 
         then
	$P(i,j_2 \onm b_2, j_2)$ and 
	$P(i,j_1 \onm b_1, j_1)$ have no common
	edges.
  \vglue -0.2 cm 
\item{(v)} Suppose
	that $v(i,j_1) \in \Omega_{j_1 \onm b_1}$, and 
	$R(j_1) \sqcap R(j_1 \onp a_1) =
	R(j_1\onm b_1) \sqcap R(j_1 \onp a_1) =
	\emptyset$.  Then
	$|P(i,j_1,j_1\onp a_1) \sqcap
	R(j_1\onm b_1)|$ is either zero or 
	greater than one.  \vglue -0.2 cm 
\item{(vi)} Suppose
	that $R(j_1)\sqcap R(j_1\onp a_1) = 
	      R(j_1\onm b_1)$ $\sqcap $ 
	$R(j_1\onp a_1) =	\emptyset$.
	Then 
	$|P(i,j_1\onm b_1,$ $j_1) \sqcap
	R(j_1\onp a_1)|$ is either zero or greater than one.

\noindent{\bitcer Proof.}
First we note that (v) and (vi) are straightforward 
consequences of the Jordan Curve Theorem.  
The arguments in the proofs of (i), (ii), (iii) and (iv)
are quite similar to each other.  For the sake of
brevity, we only
prove (iii) and (iv).

\noindent{\it Proof of (iii)}. 
It clearly suffices to show that the closed paths
$\overline{P}(i,j_2\onm b_2,j_2),
\overline{P}(i,j_1,j_1\onp $ $a_1)$ have no vertex in
common.  Suppose that 
$v(i,k)$ is in both closed paths.  Then 
$k=j_2\onm x=j_1\onp y$, where $x\in\{0,\ldots,b_2\}$, 
$y\in\{0,\ldots,a_1\}$. Thus  
$0=j_1\onm j_2\onp x\onp y=b_1\onp a_2\onp x\onp y$, and so
$b_1+y+a_2+x=kn$. Since $b_1, a_2\geq 1$, then $b_1+y+a_2+x\geq n$.
But $b_1+y\le b_1 + a_1 \le b(j_1) + a(j_1)  
 <n/2$ and $a_2+x \le a_2 + b_2 \le  a(j_2) + b(j_2)
<n/2$, since $\dd$ is robust. 
Thus we obtain a contradiction.

\noindent{\it Proof of (iv)}. 
It clearly suffices to show that 
	$\overline{P}(i,j_2 \onm b_2, j_2)$ and 
	$\overline{P}(i,j_1 \onm b_1, j_1)$ 
have no vertex in
common other than (possibly)  
$v(i,j_2)$, which might be equal to $v(i,j_1\onm b_1)$.  
Suppose that 
$v(i,k)$ is in both closed paths, for some
$k \notin\{j_2, j_1\onm b_1\}$.  Then 
$k=j_1\onm x =j_2\onm y$, where
$x\in \{0,1,\ldots,b_1-1\}$, $y\in\{1,\ldots,b_2\}$.
Let $z = j_1 \onm j_2$. 
Since $z\leq n/2$ and $y\le b_2 \le b(j_2) <n/2$
(since $\dd$ is robust), then $x=z\onp y=z+y$.
Therefore $0 \le z = x - y < b_1$.  
On the other hand, 
$j_1\onm b_1\onm j_2=kn+z-b_1$. 
Since $-n/2<z-b_1<0$, then $k=1$.
Thus
$j_1\onm b_1\onm j_2>n/2$, contradicting the
assumption that  $j_2\preceq j_1\onm b_1$. 
\hal\ms

% ***************** S E C T I O N     6 **********************

% ***********  Definition of the \ii_j             ***********

\vskip 0.5 cm

\centerline{\scacerotro 5. THE SET OF CROSSINGS ASSOCIATED TO A RED CYCLE}

\vskip 0.35 cm

% ***** I n t r o d u c t o r y      P a r a g r a p h *****

{\leftskip = 30 pt \rightskip = 30 pt
	\it 
	\noindent Throughout this section, 
	$\dd$ is a fixed robust good drawing of $\cmn$.  
	Thus, all the observations, conventions, and results 
	from Section 4  apply.
\par}

\ms

	The aim in this section is to define, for
	each red cycle $R(j)$, a set $\ii(j)$ of 
	crossings associated to $R(j)$.   

	First we define, for each red cycle $R(j)$, 
	a partition 
	$\{ \cc_j^+, \cc_j^-, \tt_j^+, \tt_j^-, \tt_j^0 \}$
	of $Z_m$, according to the following rules:\ms

\item{(i)} $i \in \cc_j^+$ iff for some
	neighborhood $N$ of $v(i,j)$,  
	$N \cap  bl(i,j\onp 1) \subseteq \Phi_j$.

\item{(ii)} $i \in \cc_j^-$ iff for some
	neighborhood $N$ of $v(i,j)$,  
	$N \cap bl(i,j\onp 1) \subseteq \Omega_j$ and
	$N \cap bl(i,j) \subseteq \Phi_j$.

\item{(iii)} $i \in \tt_j^+$ iff for some
	neighborhood $N$ of $v(i,j)$,  
	$N \cap \bigl(P(i,j \onm 1,j
	\onp 1)\setminus\{v(i,j)\}\bigr) 
	\subseteq \Omega_j$, and
	$R(j) \sqcap P(i,j,j\onp \bb)\neq\emptyset$.

\item{(iv)} $i \in \tt_j^-$ iff $i\notin\tt_j^+$ and, 
	for some
	neighborhood $N$ of $v(i,j)$,  
	$N \cap \bigl(P(i,j \onm 1,j
	\onp 1)\setminus\{v(i,j)\}\bigr) 
	\subseteq \Omega_j$, and
	$R(j) \sqcap P(i,j\onm\bb,j)\neq\emptyset$.

\item{(v)} $i \in \tt_j^0$ iff  for some
	neighborhood $N$ of $v(i,j)$,  
	$N \cap 
	 \bigl(P(i,j\onm 1, 
	j\onp 1)\setminus\{v(i,j)\}\bigr) 
	\subseteq \Omega_j$, and
	$R(j) \sqcap P(i,j\onm \bb,j\onp \bb) = 
	\emptyset$. 

\ms

	Each $\tt_j^0$ is, in turn, partitioned into 
	subsets $\tt(\beta,j)$.  To define 
	these sets we need some more 
	notation. For each $i\in Z_m, j\in Z_n$, let
	$\ob(i,j) = \min\{b\in Z_n \ | \ 
	v(i,j) \in \Omega_{j\onm  b}\}$.  Clearly,
	$\ob(i,j) \le b(j)$ for every $i$ and $j$.   
	For each  $j\in Z_n$, let 
	$\ss(j) = \{ \beta,\  0 < \beta \le b(j) \ | \ 
		    \beta = \ob(i,j) \mybox{for some} 
			i \in Z_m\}
	$.  
	For each $\beta\in \ss(j)$, define
	$\tt(\beta,j)  =
	\{
		i \in \tt_j^0 \ |\
		\ob(i,j) = \beta
	\}$.
	
	\ms

\proclaim Proposition 7.
	The following statements hold for each $j \in Z_n$.\vglue -0.1 cm 
\item{(i)}
	If $i\in\cc_j^+$, then
	$R(j)\sqcap P(i,j,j\onp  a(j)) \neq \emptyset$.\vglue - 0.2 cm 
\item{(ii)}
	If $i\in\cc_j^-$, then
	$R(j)\sqcap P(i,j\onm  b(j),j) \neq\emptyset$. \vglue - 0.2 cm 
\item{(iii)}
	If $i\in\tt_j^+$, then
	$R(j)\sqcap P(i,j,j\onp  \bb) \neq\emptyset$.\vglue - 0.14 cm 
\item{(iv)}
	If $i\in\tt_j^-$, then
	$R(j)\sqcap P(i,j\onm\bb,j) \neq\emptyset$.\vglue - 0.14 cm 
\item{(v)}
	If
	$\tt(\beta,j) \cap \tt(\beta',j')\neq \emptyset$,
	and
	$j \prec  j'$, then
	$j \preceq j' \onm   \beta'$.

\noindent{\bitcer Proof.}
	Suppose $i\in\cc_j^+$. 
	Then
	$P(i,j,j\onp  a(j)) \cap \Phi_j\neq\emptyset$. 
	On the other hand, $R(j\onp  a(j)) \subseteq \Omega_j$, 
	and so 
	$P(i,j,j\onp  a(j)) \cap \Omega_j \neq\emptyset$.  
	Since $P(i,j,j\onp  a(j))$ contains no vertex in $R(j)$,
	it follows that
	$P(i,j,j\onp  a(j))$ must cross $R(j)$.  Thus
	(i) follows.  Statement
	(ii) is proved similarly,
	and (iii) and (iv) follow  from the definitions
	of $\tt_j^+$ and $\tt_j^-$, respectively.

	Assume that 
	$i\in \tt(\beta,j) \cap \tt(\beta',j')$,
	and
	$j \prec  j'$.  
	Suppose that $j \not\preceq j'\onm  \beta'$.  Then  
	$j' \onm   \beta' \prec j$, and so
	$j = j' \onm k$ for some $k$,
	$0 < k < \beta'$.  By the definition of
	$\ob(i,j')$, $v(i,j')\in \Phi_{j}$.  But then
	$P(i,j,j')\sqcap R(j)\neq\emptyset$, since
	$i\in\tt_j^0$.
	Each edge in $P(i,j,j')$ is in $ P(i,j,j\onp  \bb)$,
	since 
	$k < \beta' \le b(j) \le \bb$.  Thus
	$P(i,j,j\onp  \bb)\sqcap R(j)\neq\emptyset$, 
	contradicting the assumption that
	$i\in \tt(\beta,j)\subseteq \tt^0_j$. 
	\hal\ms

	For each $j\in Z_n$, and each $\beta\in\ss(j)$,
	let
	$\xx(\beta,j)$ denote the set of crossings
	of the following types:\ms

\item{(A)} all the crossings between $R(j\onm  \beta)$ and $R(j)$;\vglue -0.15 cm 
\item{(B)} if $i\in\tt(\beta,j)$ and
	$R(j\onm  \beta)\sqcap P(i,j,j\onp  a(j))\neq\emptyset$,
	the last $(i,j\onm  \beta)$--crossing 
	from $v(i,j\onp  a(j))$;\vglue -0.15 cm 
\item{(C)} if $i\in\tt(\beta,j)$ and
	$R(j\onp  a(j))\sqcap P(i,j\onm  \beta,j)\neq\emptyset$,
	the first $(i,j\onp  a(j))$--crossing 
	from $v(i,j\onm  \beta)$;\vglue -0.15 cm 
\item{(D)} if $i, i' \in \tt(\beta,j)$, $i\neq i'$, 
	every crossing between
	$P(i,j\onm  \beta,j)$ and 
	$P(i',j\onm  \beta,j)$, and
	every crossing between
	$P(i,j\onm  \beta,j)$ and  
	$P(i',j,j\onp  a(j))$. 
\ms

	For each $j\in Z_n$, let $\yy(j)$ denote the collection of crossings
	of the following types:
\ms
\item{(I)}  for each $i\in \cc_j^+\cup\tt_j^+$, 
	the first $(i,j)$--crossing  
	from $v(i,j)$; and \vglue -0.2 cm 

\item{(II)} for each $i\in \cc_j^-\cup \tt_j^-$, the last
	$(i,j)$--crossing from $v(i,j)$.
\ms

	We are now ready to define the set $\ii(j)$ 
	of crossings associated to each red cycle $R(j)$:
	$$
	\ii(j)
	=
	\yy(j) \cup \biggl(\  \bigcup_{\beta \in \ss(j)}
				\xx(\beta,j) 
		    \biggr). \eqno{(1)}
	$$

	In the next section we show that if $j\neq k$, then
	$\ii(j)\cap\ii(k)=\emptyset$.

% ***************** S E C T I O N     7 **********************

% ***********  The sets \ii_j are well-defined     ***********

\vskip 0.5 cm
\ms

\centerline{\scacerotro 6. NO CROSSING IS ASSOCIATED TO MORE
	THAN ONE RED CYCLE}

\vskip 0.35 cm

% ***** I n t r o d u c t o r y      P a r a g r a p h *****

{\leftskip = 30 pt \rightskip = 30 pt
	\it 
	\noindent Throughout this section, 
	$\dd$ is a fixed robust good drawing of $\cmn$.  
	Thus, all the observations, conventions, and results 
	from Sections 4 and 5  apply.
\par}

\ms

	The main result in this section is the following.

\proclaim Lemma 8. If	$j \neq k$, then
	$\ii(j) \cap \ii(k) = \emptyset$.  That is,
	no crossing in $\dd$ is associated to more than one red
	cycle.

\noindent{\bitcer Proof.} This is an immediate consequence of
	Proposition 9 below. \hal\ms

\proclaim Proposition 9.
	Let $j, k\in Z_n$,  $\beta\in \ss(j)$, and 
	 $\beta'\in\ss(k)$.  Then:
\item{(a)} If $j\neq k$, then 
	$\yy(j)\cap\yy(k)=\emptyset$. \vglue -0.2 cm 
\item{(b)} 
	$\yy(j)\cap \xx(\beta',k)=\emptyset$. 
	\vglue -0.2 cm 
\item{(c)} If $j\neq k$ or $\beta\neq \beta'$, then
	$\xx(\beta,j)\cap\xx(\beta',k)=\emptyset$.

\noindent{\bitcer Proof of (a).} 
	Suppose  $j\neq k$.
	Each crossing in $\yy(j)$ (respectively
	$\yy(k)$) is a bichromatic 
	crossing whose red  edge involved is in $R(j)$
	(respectively $R(k)$). 	Since $j \neq k$, (a) 
	follows.

\noindent{\bitcer Proof of (b)}.  Seeking a contradiction,
	suppose that for some $j, k\in Z_n$, 
	$\beta'\in\ss(k)$, some (necessarily bichromatic,
	by the definition of $\yy(j)$) 
	crossing $x$ belongs to both 
	$\yy(j)$ and $\xx(\beta',k)$.
	Let $\ol(x)$ denote the blue edge involved in $x$. 
	Since $x$ is in $\yy(j)$, the red edge involved in $x$
	is in $R(j)$.  On the 
	other hand, a bichromatic
	crossing in $\xx(\beta',k)$ involves an
	edge in $R(j)$ only if $j$ is either 
	$k\onm  \beta'$ or $k\onp  a(k)$. Thus,
	either
	$j = k \onm   \beta'$ or $j = k\onp  a(k)$. 
	We analyze these cases separately.

\noindent{\bf Case 1.} $j = k \onm   \beta'$.  
	By the definition
	of $\xx(\beta',k)$, $\ol(x)$ is in 
	$P(i,k,k\onp  a(k))$.  Moreover, $x$ is the
	last $(i,k\onm  \beta')$--crossing (that is,
	$(i,j)$--crossing) from $v(i,k\onp  a(k))$.  
	Since $x$ is in $\yy(j)$, 
	$i$ is either in $\cc_j^-\cup\tt_j^-$ or in $\cc_j^+ \cup
	\tt_j^+$.

	Suppose that $i\in \cc_j^-\cup\tt_j^-$.  Then,
	$x$ occurs between $R(j)$ and
	$P(i,j\onm  \bb,j)$.  Thus $\ol(x)$ is in both
	$P(i,j\onm  \bb,j)$ and 
	$P(i,k,k\onp  a(k))$.  But this is
	impossible, since these open blue paths have no
	edges in common, by (i) in Proposition 6.

	Suppose now that $i\in \cc_j^+\cup\tt_j^+$.  
	Thus $x$ is the first 
	$(i,j)$--crossing from $v(i,j)$.  
	On the other hand, 
	 $x$ is the last $(i,j)$--crossing 
	from $v(i,k \onp   
	a(k))$,
	and  $\ol(x)$ 
	is in $P(i,k,k\onp$ $  a(k))$.  
	It follows from the definitions of $\tt(\beta',k)$
	and $a(k)$,  and (v) in Proposition 6, 
	that 
	$P(i,k,k\onp  a(k))$ crosses $R(j)$ at
	least twice (since they cross at least once).
	Thus we obtain a contradiction:
	$x$ cannot be at the same time the
	first $(i,j)$--crossing 
	from $v(i,j)$ and the last $(i,j)$--crossing
	from $v(i,k\onp  a(k))$, since these
	crossings are different.

\noindent{\bf Case 2.} $j = k \onp   a(k)$.  
	By the definition
	of $\xx(\beta',k)$, $\ol(x)$ is in 
	$P(i,k\onm  \beta',k)$.  
	Moreover, $x$ is the
	first $(i,k\onp  a(k))$--crossing (that is,
	$(i,j)$--crossing) from $v(i,k\onm  \beta')$.  
	Since $x$ is in $\yy(j)$, 
	$i$ is either in $\cc_j^-\cup\tt_j^-$ or in 
	$\cc_j^+\cup\tt_j^+$.  

	Suppose that $i\in\cc_j^-\cup\tt_j^-$.  Then $x$ is
	the last $(i,j)$--crossing from
	$v(i,j)$.  
	On the other hand, 
	 $x$ is the first $(i,j)$--crossing 
	from $v(i,k \onm\beta')$,
	and  $\ol(x)$ 
	is in $P(i,k\onm\beta',k)$.  
	It follows from the definitions of $\tt(\beta',k)$
	and $a(k)$, and (vi) in Proposition 6, 
	that 
	$P(i,k\onm\beta',k)$ crosses $R(j)$ at
	least twice (since they cross at least once).  
	Thus we obtain a contradiction:
	$x$ cannot be at the same time the 
	last $(i,j)$--crossing 
	from $v(i,j)$ and the first $(i,j)$--crossing
	from $v(i,k\onm\beta')$, since these
	crossings are different.

	Suppose that $i\in \cc_j^+\cup\tt_j^+$.  Then,
	$x$ occurs between $R(j)$ and
	$P(i,j,j\onp\bb)$.  Thus $\ol(x)$ is in both
	$P(i,j,j\onp\bb)$ and 
	$P(i,k\onm\beta',k)$.  
	But this is
	impossible, since these open blue paths have no
	edges in common, by (ii) in Proposition 6.

\noindent{\bitcer Proof of (c)}.  We derive a
	contradiction from the assumption that the
	following hold:
	(i) either $j\neq k$ or $\beta\neq \beta'$; and 
	(ii) there is a crossing $x$ in both
	$\xx(\beta,j)$ and $\xx(\beta',k)$.

	It follows from the definitions of $\xx(\beta,j)$
	and $\xx(\beta',k)$ that if $j = k$ and
	$\beta \neq \beta'$, then no crossing can belong to both
	$\xx(\beta,j)$ and $\xx(\beta',k)$.  
	Thus we can assume without loss of generality that
	$j\prec k$.  If $x$ involves only red edges, then
	$x\in R(j\onm  \beta)\sqcap R(j)$ 
	and $x\in R(k\onm  \beta')\sqcap R(k)$.  This clearly cannot
	happen, since there are at least three different
	cycles in $\{R(j\onm  \beta),R(j), R(k\onm  \beta'),R(k)\}$.  
	Thus $x$ involves at least an edge from a
	blue cycle $B(i)$, such that 
	$i\in \tt(\beta,j)$ and $i \in \tt(\beta',k)$. 
	By Statement (v) in Proposition 7,
	$j \preceq k \onm   \beta'$. 
	The crossing $x$ involves either one blue edge and one red
	edge or two blue edges.  We analyze these cases
	separately.

\noindent{\bf Case 1.} {\it $x$ involves one blue 
	edge $\ol(x)$ and one red
	edge $\nnor(x)$.}  Each bichromatic crossing in 
	$\xx(\beta,j)$ (respectively $\xx(\beta',k)$)
	involves a red edge in either $R(j\onm  \beta)$ or 
	$R(j\onp  a(j))$ (respectively $R(k\onm  \beta')$ or
	$R(k\onp  a(k))$).  Therefore, since 
	by assumption $x$ is in both $\xx(\beta,j)$ 
	and $ \xx(\beta',k)$,
	it follows that either 
	(a) $j\onm  \beta = k\onm  \beta'$; or
	(b) $j\onm  \beta = k\onp  a(k)$; or
	(c) $j\onp  a(j) = k\onm  \beta'$; or
	(d) $j\onp  a(j) = k\onp  a(k)$.  
	Since $j\preceq k\onm\beta'$ and $\beta<n/2$,
	it follows that $j\onm  \beta \neq 
	k\onm  \beta'$.  Therefore (a) cannot hold.
	Thus we analyze (b), (c), and (d).

	Suppose that
	$j\onm  \beta = k\onp  a(k)$.  
	It follows from
	the definitions of $\xx(\beta,j)$ and $\xx(\beta',k)$
	that $\ol(x)$ is in both
	$P(i,j,j\onp  a(j))$ and $P(i,k\onm  \beta',k)$.
	This contradicts (iii) in Proposition~6, and so
	(b) cannot hold.
	A similar argument shows that (c) cannot hold either.

	Finally, suppose that
	$j\onp  a(j) = k \onp   a(k)$.  
	It follows from
	the definitions of $\xx(\beta,j)$ and $\xx(\beta',k)$
	that $\ol(x)$ is in both
	$P(i,j\onm  \beta,j)$ and $P(i,k\onm  \beta',k)$.
	This contradicts (iv) in Proposition 6. 
        Thus (d) cannot hold.

\noindent{\bf Case 2.} {\it $x$ involves two blue edges}.
	By the definition of $\xx(\beta,j)$, $x$
	occurs between edges in different
	blue cycles $B(i)$ and $B(i')$. 
	We can assume without loss of generality that 
	$x$
	involves an edge in $P(i,j\onm  \beta,j)$, and an
	edge in either $P(i',j\onm  \beta,j)$ or
	$P(i',j,j\onp  a(j))$.  On the other hand,
	since $x$ is in $\xx(\beta',k)$, 
	$x$ 
	involves either 
	(i) an edge in $P(i,k\onm  \beta',k)$ and an
	edge in $P(i',k\onm  \beta',k)$; or
	(ii) an edge in $P(i,k\onm  \beta',k)$ and an
	edge in $P(i',k,k\onp  a(k))$; or
	(iii) an edge in $P(i',k\onm  \beta',k)$ and an
	edge in $P(i,k,k\onp  a(k))$.   
	Now (i)
	and (ii) cannot hold, since by (iv)
	in Proposition 6 
	$P(i,j\onm  \beta,j)$ has no edge in common
	with $P(i,k\onm  \beta',k)$.   On the other hand,
	it is easily checked that
	(iii) holds only if 
	(I) 
        $P(i',k\onm\beta', k)$ and $P(i',j,j\onp a(j))$ 
	have a common edge and (II) 
        $P(i,k,k\onp a(k))$ and $ P(i,j\onm\beta,j)$ have
	a common edge.  
	But (I) and (II) hold simultaneously 
	only if $\beta+\beta'+a(j)+a(k)>n$.  But
	this is impossible, since the assumption that
	$\dd$ is robust implies that 
	$\beta + a(j)\le b(j) + a(j) < n/2$ 
	and $\beta' + a(k) \le b(k) + a(k) < n/2$. \hal

\proclaim Corollary 10. Each of the unions on the 
	right hand side of Eq.~(1) is
	a disjoint union. \hal

% ***************** S E C T I O N     8 **********************

% ***  Each red cycle has associated m-2 or more crossings ***

\vskip 0.5 cm

\centerline{\scacerotro 7. AT LEAST {\bituno m} -- 2 CROSSINGS ARE
	ASSOCIATED TO EACH RED CYCLE}
\vskip 0.35 cm

% ***** I n t r o d u c t o r y      P a r a g r a p h *****

{\leftskip = 30 pt \rightskip = 30 pt
	\it 
	\noindent Throughout this section, 
	$\dd$ is a fixed robust good drawing of $\cmn$.  
	Thus, all the observations, conventions, and results 
	from Sections 4, 5,  and 6  apply.
\par}

\ms

\ms

	The purpose of  this section is to prove the following. 	

\proclaim Lemma 11. For each $j \in Z_n$,
	$|\ii(j)|\ge m - 2$.   In other words, there are at least
	$m - 2$ crossings associated to each red cycle.

\noindent{\bitcer Proof.}
	Suppose that $b(j)\in\ss(j)$.  Then,
	by Corollary 10, 
	$|\ii(j)| =  
	|\yy(j)|  +
	\bigl(\sum_{\beta\in\ss(j),\beta\neq b(j)}$ $
		|\xx(\beta,j)|\bigr) + 
		|\xx(b(j),j)|$.  Applying (i), (ii), 
	and (iii) in Proposition 12
	below,
	$|\ii(j)| \ge
	|\cc_j^+| + |\cc_j^-| + |\tt_j^+| + |\tt_j^-| +
	\bigl(\sum_{\beta\in\ss(j),\beta\neq b(j)}
		|\tt(\beta,j)|\bigr) + 
		|\tt(b(j),j)| - 2.$	
	Since $Z_m$ is the disjoint union of
	$\cc_j^+,   \cc_j^-,   \tt_j^+, \tt_j^-$, 
	and the sets $\tt(\beta,j)$ (for all
	$\beta\in \ss(j)$), it follows that
	$|\ii(j)| \ge
	m - 2$, as required. 

	Now suppose $b(j)\notin\ss(j)$.  Then,
	by Corollary 10, 
	$|\ii(j)| =  
	|\yy(j)|  +
	\bigl(\sum_{\beta\in\ss(j),\beta\neq b(j)}$
		$|\xx(\beta,j)|\bigr)$.  Applying 
	(i) and (ii) in 
	Proposition 12
	((iii) does not apply, since
	$b(j)\notin\ss(j)$),
	$|\ii(j)| \ge
	|\cc_j^+| + |\cc_j^-| + |\tt_j^+| + |\tt_j^-| +
	\bigl(\sum_{\beta\in\ss(j),\beta\neq b(j)}
		|\tt(\beta,j)|\bigr)$.
	Since $Z_m$ is the disjoint union of
	$\cc_j^+,   \cc_j^-,   \tt_j^+, \tt_j^-$, 
	and the sets $\tt(\beta,j)$ (for all
	$\beta\in \ss(j)$), we obtain
	$|\ii(j)| \ge
	m$. \hal\ms

\proclaim Proposition 12. For each $j\in Z_n$, 
	the following statements hold.
		\item{(i)} 
	$|\yy(j)| \ge |\cc_j^+| + |\cc_j^-| +
			|\tt_j^+| + |\tt_j^-|$.
			 \vglue -0.15 cm 
	\item{(ii)}	For  each
	$\beta\in \ss(j)$,$\beta\neq b(j)$, 
	$|\xx(\beta,j)| \ge |\tt(\beta,j)|$.	\vglue -0.15 cm 
	\item{(iii)} 	
	If $b(j)\in \ss(j)$, then
	$|\xx(b(j),j)| \ge |\tt(b(j),j)| - 2$.

\noindent{\it Proof of (i)}.
	If $i\in\cc_j^+$, then by 
	Proposition 7 
	$P(i,j,j\onp  a(j)) \sqcap
	R(j)\neq\emptyset$, and by the definition 
	of $\yy(j)$, one of these crossings is in $\yy(j)$.  
	If $i\in \cc_j^-$, then by 
	Proposition 7
	$P(i,j\onm  b(j),j) \sqcap
	R(j)\neq\emptyset$, and by the definition
	of $\yy(j)$, one of these crossings is in $\yy(j)$. 
	If $i\in\tt_j^+$, then by
	Proposition 7
	$P(i,j,j\onp  \bb)\sqcap  
	R(j)\neq\emptyset$, and by the definition
	of $\yy(j)$, one of these crossings is in $\yy(j)$. 
	If $i\in\tt_j^-$, then by
	Proposition 7
	$P(i,j\onm \bb,j)\sqcap
	R(j)\neq\emptyset$, and by the definition
	of $\yy(j)$, one of these crossings is in $\yy(j)$. 
	Since $\cc_j^+, \cc_j^-, \tt_j^+$, and $\tt_j^-$ are
	pairwise disjoint, (i) follows.  \ms

\noindent{\it Proof of (ii)}.
	Let $\beta\in\ss(j), \beta\neq b(j)$. 
	Let $\aa = \{P(i,j\onm  \beta,j\onp  a(j)\ | \
	i \in \tt(\beta,j)$).
	Since $\beta < b(j)$, $R(j\onm  \beta) \sqcap R(j) = k > 0$. 
	Thus it follows from the definition of 
	$\tt(\beta,j)$ that 
	$(\{R(j \onm   \beta),R(j),R(j\onp  a(j))\}, \aa\}$ 
	is a $k$--intersecting 
	$(3,|\tt(\beta,j)|)$--configuration.
	Note that the crossings in 
	$(\{R(j\onm  \beta),R(j),R(j\onp  a(j))\}, \aa)$ that
	are in $\xx(\beta,j)$ are (a) the good
	crossings; (b) one crossing for each initial
	arc that crosses $R(j\onp  a(j))$; 
	(c) one crossing for each final arc that
	crosses $R(j\onm  \beta)$; and
	(d) the $k$ crossings between $R(j\onm  \beta)$ and $R(j)$. 
	By Corollary 5, 
	there are at least $|\aa| = 
	|\tt(\beta,j)|$ such crossings. \ms

\noindent{\it Proof of (iii)}.
	Suppose that $b(j)\in\ss(j)$. 
	Let $\aa = \{P(i,j\onm  b(j),j\onp  a(j)\ | \
	i \in \tt(b(j),j)$).
	By the definition of $b(j)$, 
	$R(j\onm  b(j)) \sqcap R(j) = 0$. 
	Thus it follows from the definition of 
	$\tt(b(j),j)$ that 
	$(\{R(j \onm   b(j)),R(j),R(j\onp  a(j))\}, \aa\}$ 
	is a $0$--intersecting 
	$(3,|\tt(b(j),j)|)$--configuration.
	Note that the crossings in 
	$(\{R(j\onm b(j)),R(j),R(j\onp  a(j))\},$ $\aa),$ 
	that are in $\xx(b(j),j)$ are (a) the good
	crossings; (b) one crossing for each initial
	arc that crosses $R(j\onp  a(j))$; and
	(c) one crossing for each final arc that
	crosses $R(j\onm  b(j))$.
	By Corollary 3, 
	there are at least 
	$|\aa| = |\tt(b(j),j)| - 2$ 
	such crossings. 
	 \hal\ms

% ***************** S E C T I O N     9 **********************

% *******  Proofs of Theorem 1 and of the Main Theorem  ******

\vskip 0.5 cm

\centerline{\scacerotro 8. PROOFS OF THEOREM 1 AND THE MAIN THEOREM}

\vskip 0.35 cm

% ***** I n t r o d u c t o r y      P a r a g r a p h *****

	First we show that
	if $n$ is sufficiently large compared to $m$, 
	then every drawing of $\cmn$  either is robust
	or has a red cycle with at least $m$ crossings.

\proclaim Proposition 13.  Let $m, n$ be such that
	$n \ge (m + 3)^2/2 + 1$, $m \ge 3$.  Let $\dd$ be
	a drawing of $\cmn$.  Then either $\dd$ is robust
	or there is a red cycle with at least $m$
	crossings in $\dd$.  

\noindent{\bitcer Proof.} 
	Let $m, n$ satisfy the inequalities in the
	statement of the proposition.  Let $\dd$ be
	a drawing of $\cmn$, and suppose that no red
	cycle has $m$ or more crossings in $\dd$. 
	We will show that then $\dd$
	is robust.  Let $\rr$ denote the set of
	all red cycles. 

	Since two red cycles that cross do so
	in at least two points, it follows that each
	red cycle crosses at most $(m-1)/2$ other
	red cycles in $\dd$.   Since $n > (m-1)/2 + 1$, 
	for each $R(j)$ there is a cycle
	in $\rr\setminus\{R(j)\}$
	that does not cross 
	$R(j)$.  This shows that $b(\dd,j)$ is defined 
	for every $j\in Z_n$.   Moreover,
	$b(\dd,j) \le (m-1)/2 + 1 = (m+1)/2$ 
	for every $j\in Z_n$.  

	For each $j\in Z_n$, let 
	$\rr(j) = \{R(j\onm b(\dd,j)),R(j\onm (b(\dd,j) - 1)),
	\ldots,R(j)\}$.  Thus,
	$|\rr(j)| \le (m+3)/2$.
	We now show that $a(\dd,j)$ exists and is
	at most $(m+1)(m+3)/4 + 1$ for every
	$j\in Z_n$.  
	Let $j\in Z_n$ be fixed.  Since every red cycle
	crosses at most $(m-1)/2$ other red cycles,
	then the collection of red cycles
	that either are in $\rr(j)$ or cross a
	cycle in $\rr(j)$ has size at most
	$((m+3)/2) + ((m+3)/2)((m-1)/2)
	= (m+1)(m+3)/4$.  
	Since
	$n >(m+1)(m+3)/4 + 1$,
	it follows that there is a red cycle
	not in $\rr(j)$ that crosses no cycle
	in $\rr(j)$. Moreover,
	it follows that $a(\dd,j) \le 
	(m+1)(m+3)/4 + 1$.

	Finally, we note that 
	$b(\dd,j) + a(\dd,j) 
	\le (m+1)/2 + (m+1)(m+3)/4  + 1
	= (m+3)^2/4 < 
	n / 2$.  Thus $\dd$ is robust. 
\hal\ms

\noindent{\bitcer Proof of Theorem 1.}
	It follows immediately from the definition of $\ii(j)$, 
	Lemma 8, and Lemma 11. \hal\ms

\noindent{\bitcer Proof of Main Theorem.}
	It is easy to exhibit drawings of $\cmn$ with
	exactly $(m-2)n$ crossings, for all $m,n $
	such that $n\ge m \ge 3$.  Thus we need to
	show that $\cr(\cmn) \ge (m-2)n$, for all
	$m,n$ such that $n \ge 
	(m/2)( (m+3)^2/2 + 1)$, 
	$m \ge 3$. 
	
	Let $m \ge 3$ be fixed. 
	Let $n_0 = (m+3)^2/2 + 1$.
	For every 
	$n \ge (m/2)( (m+3)^2/2 + 1)$, 
	$\min\{(m-2)n, m(n- n_0 )\} = (m-2)n$.
	Therefore
 	it suffices to show
	that if $n \ge n_0$, then
	$\cr(\cmn) \ge \min\{(m-2)n, 
	m(n- n_0 )\}$.
	We prove this by induction on $n$. 

	The base case
	is $n = n_0$, for which there is nothing to
	prove. Suppose that the statement holds for
	$n = k-1 \ge n_0$, and consider a drawing 
	$\dd$ of $C_m\times C_{k}$.  If $\dd$ is
	robust, then we are done, since
	by Theorem 1 $\dd$ has at least
	$(m-2)k$ crossings.  Thus we  
	assume $\dd$ is not robust.   Since 
	$k > n_0$, it follows from Proposition 13
	that there is a red cycle $R$ with $m$ or
	more crossings. The drawing $\dd'$ that results
	by removing $R$ from $\dd$ has, by the induction
	hypothesis, at least
	$\min\{ (m-2)(k-1), m(k-1 - n_0)$ 
	crossings, and  so $\dd$ has at least 
	$\min\{ (m-2)(k-1) + m, m(k-1-n_0) + m\} =
	\min\{ (mk-2(k-1),  m(k-n_0)\}$ crossings.
	Since $mk -2(k-1) > (m-2)k$, then $\dd$ has
	at least 
	$\min\{ (m-2)k, m(k-n_0)\}$ crossings, as
	required. \hal\ms

% ***************** S E C T I O N     9 **********************

% *******   C o n c l u d i n g     R e m a r k s        ******

\vskip 0.5 cm

\centerline{\scacerotro 9. CONCLUDING REMARKS}

\vskip 0.35 cm

% ***** I n t r o d u c t o r y      P a r a g r a p h *****

	As we mentioned in Section 1, 
	although we do not use
	the notion of an arrangement explicitly, our proof of the
	Main Theorem is largely based on 
	the theory 
	of arrangements introduced by Adamsson [1] 
	and further developed by 
	Adamsson and Richter [2].  

	An $(m,n)$--circular arrangement consists
	of two collections $\bb, \rr$ of (blue and
	red, respectively) closed curves.
	The red curves are cyclically ordered, and
	each blue curve intersects the red curves
	in the given cyclic order.  Clearly, each
	drawing of
	$\cmn$ yields an $(m,n)$--circular 
	arrangement, if we regard each vertex as
	an intersection between a red curve and a 
	blue curve.
	Moreover, lower bounds on the number
	of intersections of $(m,n)$--circular arrangements
	imply lower bounds for $\cr(\cmn)$.  Using this
	approach, Adamsson proved that large classes of
	drawings of $\cmn$ have at least $(m-2)n$ crossings.
	He also used this approach to show that
	$\cr(C_7\times C_n) = 5n$, as conjectured.

	A nice aspect of the Adamsson and Richter approach 
	to the HKS--conjecture is that,
	although it draws from and generalizes ideas
	introduced in previous work, it is virtually
	self--contained.  A similar observation holds for 
	this paper, where no results from previous
	work on the HKS--conjecture are used to prove 
	the Main Theorem.

	In the proof of the Main Theorem we used as a base case
	of the induction the (obviously true) inequality
	$\cr(C_m\times C_{n_0}) \ge 0$.  
	It is natural to ask whether the statement of the
	Main Theorem is substantially improved if instead we
	use a nontrivial bound for $\cr(C_m\times C_{n_0})$. 
	The best general lower bound known for 
	the crossing number of $\cmn$ (for $n\ge m \ge 3$)
	is  $\cr(C_m\times C_n) \ge (1/2)(m-2)n$ [8].  
	Using this bound, we obtain the following
	slightly improved version of the Main Theorem.

\proclaim Main Theorem [Improved version].
	Let $m, n$ be integers such that
	$n \ge 
	(m/4 + 1/2)( (m+3)^2/2 + 1)$, 	
	$m \ge 3$. Then
	$\cr(\cmn) = (m-2)n$.

	As we mentioned above, $(m,n)$--circular
	arrangements are more general structures than
	drawings of $\cmn$.  Thus, Adamsson's results
	actually imply lower bounds for the crossing
	numbers of families of graphs more general than
	$\cmn$.   A similar observation holds for the
	work in this paper.  Let us say that
	a $4$--regular graph $G$ is
	an $(m,n)$--{\it graph} if $G$ consists 
	of $n$ pairwise disjoint,
	cyclically ordered (red) $m$--cycles
	$\{R(0),R(1),\ldots,R(n-1)\}$, plus 
	$mn$ (blue) edges, such that 
	for each vertex $v$ in $R(j)$, one blue edge
	incident with $v$ is incident with $R(j-1)$, and
	the other blue edge incident with $v$ is
	incident with $R(j+1)$. 
	It can be checked that the techniques developed
	above yield the following more general version
	of the Main Theorem.

\proclaim Theorem.
	Let $m, n$ be integers such that
	$n \ge 
	(m/2)((m+3)^2/2 + 1)$, $m \ge 3$. 
	Let $G$ be an $(m,n)$--graph.  Then
	$\cr(G)\ge (m-2)n$.

	While the Main Theorem settles the HKS--Conjecture
	for all but finitely many values of $n$, for each
	$m$, the HKS--Conjecture remains open for
	$n <  (m/2)((m+3)^2/2 + 1)$, $m \ge 8$.  
	For
	values of $n$ sufficiently close to $m$ (more
	precisely, for $m, n$ such that $m \ge 8$, 
	$m \le n \le 5(m-1)/4$),  
	it is known that $\cr(\cmn) \ge (5/7)mn$ [13]. 
	For $n$ between $5(m-1)/4$ and
	$(m/4 + 1/2)( (m+3)^2/2 + 1)$, 	the best 
	general lower
	bound known is 
	$cr(\cmn) \ge (m-2)n/2$ [8].

% ***************** Acknowledgements **********************

\vskip 0.5 cm

\centerline{\scacerotro ACKNOWLEDGEMENTS}

\vskip 0.35 cm

	This work was done while the first author visited
	the second at IICO--UASLP.  This visit
	was partially funded by 
	the second author's CONACYT Grant J32168E.  The
	second author also acknowledges partial 
	support from FAI--UASLP.

% **********************************************************

% ***************** Acknowledgements **********************

\vskip 0.5 cm

\centerline{\scacerotro REFERENCES}

\vskip 0.35 cm

\item{[1]} J.~Adamsson, Ph.D.~Thesis.  Carleton University (2000). 

\item{[2]} J.~Adamsson and R.B.~Richter, 
Arrangements, Circular Arrangements and the Crossing Number of
$C_7\times C_n$, manuscript. 

\item{[3]} M.~Anderson, R.B.~Richter, and 
P.~Rodney, The Crossing Number of $C_6 \times C_6$,
{\it Cong.~Numerantium} {\bf 118} (1996), 97--107. 

 \item{[4]} M.~Anderson, R.B.~Richter, and
P.~Rodney, The Crossing Number of $C_7 \times C_7$, 
{\it Cong.~Numerantium} {\bf 125} (1997), 97--117.

\item{[5]} L.W.~Beineke and R.D.~Ringeisen, On
the crossing numbers of products of cycles and graphs
of order four, {\it J.~Graph Theory} {\bf 4} (1980), 145--155.

 \item{[6]} A.M.~Dean and R.B.~Richter, The Crossing
Number of $C_4 \times C_4$, {\it J.~Graph Theory} 
{\bolcer 19} (1995), 125--129.

 \item{[7]} F.~Harary, P.C.~Kainen, and 
A.J.~Schwenk,
Toroidal graphs with arbitrarily high crossing numbers, 
{\it Nanta Math.} {\bolcer 6} (1973), 58--67.

\item{[8]} H.A.~Juarez and G.~Salazar,
	Drawings of $\cmn$ with one disjoint family II,
	{\it J.~Combinat.~Theory Series B},
	to appear.

\item{[9]} M.~Kle{\v s}{\v c}, R.B.~Richter and
I.~Stobert, The crossing
number of $C_5 \times C_n$, 
{\it J.~Graph Theory} {\bolcer 22} (1996), 239--243.

\item{[10]} R.B.~Richter and G.~Salazar, The crossing number
	of $C_6\times C_n$, {\it Australasian Journal of
	Combinatorics}, to appear.

 \item{[11]} R.B.~Richter and C.~Thomassen, Intersections of 
Curve Systems and the Crossing Number of $C_5 \times C_5$, 
{\it Disc.~Comp.~Geom.} {\bolcer 13} (1995), 149--159.

 \item{[12]} R.D.~Ringeisen and L.W.~Beineke, 
The crossing number of $C_3 \times C_n$, 
{\it J.~Combinat.~Theory} {\bolcer 24} (1978), 134--136.

\item{[13]} F.~Shahrokhi, O.~S\'ykora, L.~A.~Sz\'ekely
and I.~Vr\v to,  Intersection of Curves and Crossing 
Number of $C_m \times C_n$ on Surfaces,  {\it Discrete Comput.~Geom.}
{\bf 19} 
(1998), no.~2, 237--247.

\end